\documentclass[12pt,twoside]{amsart}

\theoremstyle{definition}

\theoremstyle{remark}

\numberwithin{equation}{section}

\date{}

\begin{document}

\centerline{\bf Int. J. Contemp. Math. Sciences, Vol. 4, 2009, no.
12, 573 - 576}

\centerline{}

\centerline{}

\centerline {\Large{\bf Convexity of \u{C}eby\u{s}ev Sets}}

\centerline{}

\centerline{\Large{\bf in Hilbert Spaces}}

\centerline{}

\centerline{\bf {Hadi Haghshenas}}

\centerline{}

\centerline{Department of Mathematics, Birjand University,
Birjand, Iran}

\centerline{h$_{-}$haghshenas60@yahoo.com}

\centerline{}

\centerline{}

{\bf Abstract.} The aim of this paper is state of conditions that
ensure the convexity of a \u{C}eby\u{s}ev sets in Hilbert spaces.

\centerline{}


\section{\textbf{Introduction}} The approximation theory
is one of the important branch of functional analysis that
\u{C}eby\u{s}ev originated it in nineteenth century. But,
convexity of \u{C}eby\u{s}ev sets is one of the basic problems in
this thoery. In a finite dimensional smooth normed linear space a
\u{C}eby\u{s}ev set is convex [2]. Also, every boundedly compact
\u{C}eby\u{s}ev set in a smooth Banach space is convex [3,7] and
in a Banach space which is uniformly smooth, each approximately
compact \u{C}eby\u{s}ev set is convex [4]. In addition, in a
strongly smooth space, every \u{C}eby\u{s}ev set with continuous
metric projection is convex [5,6]. Regarding convexity of
\u{C}eby\u{s}ev sets, there are still several open problems. It is
a wellknown problem whether a \u{C}eby\u{s}ev set in a Hilbert
space must be convex. Of course, in a finite-dimensional Hilbert
space, every \u{C}eby\u{s}ev set is convex. For inverse, we know
that every closed convex set in a strictly convex reflexive Banach
space and in particular Hilbert space is \u{C}eby\u{s}ev. However,
this problem that whether every \u{C}eby\u{s}ev set in a strictly
convex reflexive Banach space is convex is still open.

\section{\textbf{Basic definitions and Preliminaries}} In this section
we collect some elementary facts which will help us to establish
our main results.\\\textbf{Definition 2.1. } Let $(X,\|.\|)$ be a
real normed linear space, $x \in X $ and $X^{*}$ be its dual
space. For a nonempty subset $K$ in X, the distance of $x$ from
$K$ is defined as $d_{K}(x)= inf \{\|x-v\| ; v\in K \}$. The $K$
is said to be a \u{C}eby\u{s}ev set if, each point in $X$ has a
unique best approximation in $K$. In other words, for every $x \in
X $, there exist a unique $v\in K$ such that $\|x-v\|= d_{K}(x)$.
(This concept was introduced by S. B. Stechkin in honour of the
founder of best approximation theory, \u{C}eby\u{s}ev). The metric
projection is given by $P_{K}(x)= \{v\in K ; \|x-v\|= d_{K}(x) \}$
which consists of the closest points in $K$ to $x$. The $P_{K}$ is
said to be continuous if, $P_{K}(x)$ is a singleton for each $x\in
X \setminus K$ and it is sequently continuous.
\\\textbf{Definition 2.2. } A norm $\|.\|$ on X is said to be Kadec if, each
weakly convergent sequence $(x_{n})_{n=1}^{\infty}$ in X with the
weak limit $x\in X$ converges in norm to $x$ whenever
$\|x_{n}\|\rightarrow\|x\|$ as $n \rightarrow
\infty$.\\\textbf{Definition 2.3. } The space $X$ is said to be
strictly convex if, $x=y$ whenever $x , y \in S(X)$ and
$\displaystyle{\frac{x+y}{2}\in S(X)}$, where $S(X)=\{x\in
X;\|x\|=1\}$.
\\\textbf{Remark 2.4. } Every Hilbert space is strictly convex. Hence the dual of each Hilbert space is strictly convex.\\Related to the notion of strict convexity, is the notion
of smoothness.\\ \textbf{Definition 2.5. }For each $x\in X$ the
element \ $x^{*}\in S(X^{*})$ satisfying
$\|x\|=\hspace{1mm}\langle x^{*},x \rangle $ is called the support
functional corresponding to $x$ and $X$ is smooth in a non-zero
$x\in X$ if, the support functional corresponding to $x$ is
unique.\\Of course, the Hahn-Banach extension theorem, ensures the
existence of at least one such support functional.\\Smoothness and
strict convexity are not quite dual properties. There are examples
of smooth spaces whose duals fail to be strictly
convex.\\\textbf{Theorem 2.6. }[1]  Each Hilbert space is
smooth.\\\textbf{Example 2.7. }The space $\Bbb{R}^{2}$ with
Euclidian norm,
$$\|(x_{1},x_{2})\|=\sqrt{x_{1}^{2}+x_{2}^{2}}$$ is a smooth and
strictly convex space.
\\\\\textbf{Theorem 2.8. }[1]
If $X$ be a reflexive and smooth space, then the dual space $X^{*}$ is strictly convex.\\
\textbf{Definition 2.9.} The space $X$ is uniformly convex if, for
every sequences $( x_{n})_{n=1}^{\infty}$ and $(
y_{n})_{n=1}^{\infty}$ we have,
$\displaystyle{\lim_{n\rightarrow\infty}\|x_{n}-y_{n}\|=0}$
whenever, $\displaystyle{\lim_{n\rightarrow\infty}\|x_{n}+y_{n}\|=2}$\\
\textbf{Theorem 2.10. }[1]  Every uniformly convex space, is
strictly convex.\\\textbf{Remark 2.11. }[1] The inverse of theorem
(2.10) is not true, necessary. For example, define a norm
$|\|.\||$ on $C[0,1]$ by
$|\|x\||^{2}=\|x\|_{\infty}^{2}+\|x\|_{2}^{2}$, where
$\|.\|_{\infty}$ and $\|.\|_{2}$, denote the norms of $C[0,1]$ and
$L^{2}[0,1]$, respectively. Then $|\|.\||$ is strictly convex but,
not uniformly convex on $C[0,1]$.\\\textbf{Theorem 2.12. }[1] Each
Hilbert space, is uniformly convex.\\\textbf{Theorem 2.13. }[1]
Every uniformly convex Banach space, is
reflexive.\\\textbf{Theorem 2.14. }[1] The norm
of every uniformly convex space, is Kadec.\\
\section{\textbf{Main Results}}In this section, we state the conditions that ensure the
convexity of a \u{C}eby\u{s}ev sets in Hilbert
spaces\\\textbf{Theorem 3.1. } Let $K$ be a weakly closed set in a
reflexive space X with Kadec norm. Then the metric projection
$P_{K}$ is continuous.
\\\textit{Proof.   }Let $x \in X\backslash K$, $v \in P_{K}(x)$
and suppose $(x_{n})_{n=1}^{\infty} \subseteq X $,
$(v_{n})_{n=1}^{\infty} \subseteq P_{K}(x_{n})$ such that
$x_{n}\rightarrow x$ in norm. It is sufficient show that
$P_{K}(x)$ is a singleton and $v_{n}\rightarrow v$ in norm. Since
that $d_{K}$ is continuous, we have :
$$d_{K}(x)= \|x-v\|\leq \|v_{n}-x\|\leq\|x_{n}-v_{n}\|+\|x_{n}-x\|=d_{K}(x_{n})+\|x_{n}-x\|\rightarrow d_{K}(x)$$
So, $\displaystyle{\lim_{n\rightarrow\infty}\|v_{n}-x\|=\|x-v\|}$
and hence $(v_{n})_{n=1}^{\infty}$is bounded. Thus
$(v_{n})_{n=1}^{\infty}$ is contained in an set $A$ such that $A$
is weakly closed and boundedly in norm. Since $X$ is reflexive,
the set $A$ is weakly compact. Hence there exists a weakly
convergent subsequence $(v_{n_{k}})_{k=1}^{\infty}$ of
$(v_{n})_{n=1}^{\infty}$ whose weak limit $v_{0}$ lies in $A$:
Such an $v_{0}$ must be in $K$. Note that the norm on a normed
space is lower semicontinuous for the weak topology. Then
$$\|x-v\|=d_{K}(x)\leq \|x-v_{0}\|\leq \displaystyle{\liminf_{k\rightarrow\infty}\|v_{n_{k}}-x\|=d_{k}(x)=\|x-v\|}$$
This implies $v_{0}=v$ and so $P_{K}(x)$ is a singleton. The
$(x-v_{n_{k}})_{k=1}^{\infty}$ is weakly converges to $x-v$ and
satisfies
$\displaystyle{\lim_{k\rightarrow\infty}\|v_{n_{k}}-x\|=\|x-v\|}$.
Since the norm on X is Kadec, the sequence
$(x-v_{n_{k}})_{k=1}^{\infty}$ is normly convergent to $(x-v)$.
Therefore, $( v_{n_{k}})_{k=1}^{\infty}$ converges to $v$ in norm
and consequently $( v_{n})_{n=1}^{\infty}$ converges to $v$ in
norm . This proves that $P_{K}$ is continuous.\\Now by the
Theorems 2.12, 2.13, 3.1, we have:\\\textbf{Corollary 3.2. } Let
$K$ be a weakly closed set in a uniformly convex Banach space X.
Then the metric projection $P_{K}$ is continuous.
\\\textbf{Theorem 3.3. }[5,6] Every \u{C}eby\u{s}ev set $K$ with
continuous metric projection $P_{K}$, in a Banach space $X$ with
strictly convex dual $X^{*}$, is convex.\\Now by Remark 2.4 and
the previous theorem, we have:\\\textbf{Corollary 3.4. } Every
\u{C}eby\u{s}ev set with continuous metric projection, in a
Hilbert space is convex.\\Now by the Theorems 3.3, 2.12, 2.8 and
Corollary 3.2, we have:\\\textbf{Theorem 3.5. } Every weakly
closed \u{C}eby\u{s}ev set in a smooth uniformly convex Banach
space, is convex.\\Finally, by the Theorems 2.6, 2.11 and the
previous theorem, we have:\\\textbf{Corollary 3.6.} Every weakly
closed \u{C}eby\u{s}ev set in a Hilbert space, is convex .\\
\section{\textbf{Acknowledgments}}
The author is indebted to his supervisor professor Amanollah
Assadi for the useful remarks while this work was in progress.\\

\centerline{}

{\bf Received: September, 2008}

\end{document}